\documentclass{article}

\usepackage{amsmath}
\usepackage{amssymb}
\usepackage{enumitem}
\usepackage{graphicx}
\usepackage{mathdots}
\usepackage{color}

\def\oM{\overline{\mathcal{M}}}
\def\cM{{\mathcal{M}}}

\def\Z{\mathbb{Z}}
\def\C{\mathbb{C}}
\def\Q{\mathbb{Q}}

\def\qed{{\hfill $\Diamond$}}
\def\b1{{\bf 1}}

\def\com{\mathbb{C}}

\usepackage{theorem}
\newtheorem{definition}{Definition}
\newtheorem{theorem}[definition]{Theorem}
\newtheorem{remark}[definition]{Remark}
\newtheorem{proposition}[definition]{Proposition}
\newtheorem{corollary}[definition]{Corollary}

\newtheorem{prop}{Proposition}
\newtheorem{defn}[prop]{Definition}
\newtheorem{problem}[prop]{Problem}
\newtheorem{rem}[prop]{Remark}

\newcommand{\MM}{\oM}

\title{Cohomological field theories with non-tautological classes}
\author{R. Pandharipande and D. Zvonkine}
\date{April 2019}

\begin{document}

\maketitle

\vspace{-20pt}

\begin{abstract} 
A method of constructing Cohomological Field Theories (CohFTs) with unit
using minimal classes on the moduli spaces of curves is developed. As
a simple consequence, CohFTs with unit are found which take values
outside of the tautological cohomology of the moduli spaces of curves.
A study of minimal classes in low genus is presented in the 
Appendix by D. Petersen.
\end{abstract}

\setcounter{tocdepth}{1} 


\setcounter{section}{-1}
\section{Introduction}

\subsection{Moduli of curves} \label{xcxc}
Let $\oM_{g,n}$ be the moduli space of Deligne-Mumford stable curves of
genus $g$ with $n$ markings \cite{DM}. 
There are natural forgetful morphisms dropping the last marking,
$$p:\oM_{g,n+1} \rightarrow \oM_{g,n}\, ,$$
and boundary morphisms
$$
q : \oM_{g-1, n+2} \to \oM_{g,n}\, ,
$$
$$
r: \oM_{g_1, n_1+1} \times \oM_{g_2, n_2+1} \to \oM_{g,n}\, , 
$$
where $n=n_1+n_2$ and $g=g_1+g_2$. The images of both $q$ and $r$
lie in the boundary $$\partial \oM_{g,n} \subset {\oM}_{g,n}\, .$$
Stability requires $2g-2+n>0$.

The cohomology and Chow groups 
of the moduli space of curves are
$$H^*({\oM}_{g,n},\mathbb{C}) \ \ \ \text{and} \ \ \  A^*({\oM}_{g,n},\mathbb{C})\, .$$
See \cite{FPan,PandSLC} for a survey of results and open questions.

\subsection{Cohomological field theories}\label{cohax}
The starting point for defining a Cohomological Field Theory \cite{KonMan} is
a triple of data $(V,\eta,\b1)$ where
\begin{enumerate}
\item[$\bullet$]
$V=V_0\oplus V_1$ is a finite dimensional 
$\mathbb{C}$-vector space{\footnote{Often
CohFTs are defined over the field $\mathbb{Q}$. However, for our
examples here, we
 will require $\C$.}} with a $\Z_2$-grading, 
\item[$\bullet$] $\eta$ is an even nondegenerate quadratic form -- 
a non-degenerate 2-form on $V$ which is symmetric on $V_0$,
skew-symmetric on $V_1$, and satisfies 
$$\eta(V_0,V_1)=\eta(V_1,V_0)=0\, ,$$
\item[$\bullet$] $\b1 \in V_0$ is a distinguished element.
\end{enumerate}
Given a 
$\mathbb{C}$-basis $\{e_i\}$ of $V$, the 
symmetric form $\eta$ can be written as a matrix
$$\eta_{jk}=\eta(e_j,e_k) \ .$$ The inverse matrix is denoted, as usual, by
$\eta^{jk}$. We will only consider bases which respect the grading of
$V$.

A {\em Cohomological Field Theory} consists of 
a system $\Omega = (\Omega_{g,n})_{2g-2+n > 0}$ of {even} tensors
$$
\Omega_{g,n} \in H^*(\oM_{g,n},\mathbb{C}) \otimes (V^*)^{\otimes n}.
$$
The tensor $\Omega_{g,n}$ associates a cohomology class in $H^*(\oM_{g,n},\mathbb{C})$ 
to vectors 
$$v_1, \ldots, v_n\in V$$ assigned to the $n$ markings. Let $\Omega_{g,n}(v_1,\ldots, v_n) $
denote the associated cohomology class in $H^*(\oM_{g,n},\mathbb{C})$. For 
vectors $v_i\in V_{\overline{v}_i}$ of pure grading
$\overline{v}_i \in \{0,1\}$,
 the {\em even} condition is
$$
\mbox{deg\,} \Omega_{g,n}(v_1,\ldots, v_n) = 
\overline{v}_1 + \cdots + \overline{v}_n \bmod 2\, ,
$$
where $\text{deg}$ denotes the cohomological degree
in  $H^*(\oM_{g,n},\mathbb{C})$.

In order to define a Cohomological Field Theory, the system 
$$\Omega = (\Omega_{g,n})_{2g-2+n > 0}$$ must satisfy the 
CohFT axioms:

\begin{enumerate}
\item[(i)] Each tensor $\Omega_{g,n}$ is $S_n$-invariant (in
the $\mathbb{Z}_2$-graded sense) for
 the natural action of the symmetric group $S_n$ on
$$H^*(\oM_{g,n},\mathbb{C}) \otimes (V^*)^{\otimes n}$$
 obtained by simultaneously 
permuting the $n$ marked points of $\oM_{g,n}$ and the $n$ factors of $V^*$. 
\item[(ii)] 
The tensor $q^*(\Omega_{g,n}) \in H^*(\oM_{g-1,n+2},\mathbb{C}) \otimes 
(V^*)^{\otimes n}$,
obtained via pull-back by the boundary morphism 
$$
q : \oM_{g-1, n+2} \to \oM_{g,n}\ ,
$$
is required to equal the contraction
of $\Omega_{g-1,n+2}$ by the bi-vector 
$$\sum_{j,k} \eta^{jk} e_j \otimes e_k$$
inserted at the two identified points: 
$$q^*(\Omega_{g,n}(v_1,\ldots,v_n)) =
\sum_{j,k} \eta^{jk}\,\Omega_{g-1,n+2}(v_1,\ldots,v_n, e_j,e_k)$$
in $H^*(\oM_{g-1,n+2},\mathbb{C})$ for all $v_i \in V$.

The tensor $r^*(\Omega_{g,n})$, obtained via pull-back 
by the boundary morphism{\footnote{Here, we
assume the $n_1$ markings of first factors are
$\{1,\ldots,n_1\}$ and the $n_2$ markings of the second
factor are $\{n_1+1,\ldots, n_1+n_2=n\}$. Such a
marking distribution can always be achieved by $S_n$-action (with  a 
possible sign change due to the grading).}}
$$
r: \oM_{g_1, n_1+1} \times \oM_{g_2, n_2+1} \to \oM_{g,n}\, , 
$$
is similarly required to equal the contraction of 
$\Omega_{g_1, n_1+1} \otimes \Omega_{g_2, n_2+1}$ by the
same bi-vector:
\begin{multline*}
r^*(\Omega_{g,n}(v_1,\ldots,v_n)) = \\
\sum_{j,k} \eta^{jk}\, \Omega_{g_1,n_1+1}(v_1,\ldots,v_{n_1}, e_j) \otimes
\Omega_{g_2,n_2+1}(e_k, v_{n_1+1}, \ldots, v_n)
\end{multline*}
in $H^*(\oM_{g_1,{n_1+1}},\mathbb{C}) \otimes 
H^*(\oM_{g_2,{n_2+1}},\mathbb{C})$ for all $v_i \in V$.

\item[(iii)] The tensor $p^*(\Omega_{g,n})$, obtained via pull-back by the forgetful map
$$p:\oM_{g,n+1} \rightarrow \oM_{g,n}\, ,$$
is required to satisfy 
$$
\Omega_{g,n+1}(v_1, \ldots,v_n , \b1) = p^*\Omega_{g,n} (v_1, \ldots,  v_n)\ 
$$
for all $v_i \in V$.
In addition, the equality 
$$\Omega_{0,3}(v_1,v_2,  \b1) = \eta(v_1,v_2)\, $$
is required for all $v_i \in V$.
\end{enumerate}

\begin{definition}\label{defcohft}
A system $\Omega= (\Omega_{g,n})_{2g-2+n>0}$ of tensors 
$$
\Omega_{g,n} \in H^*(\oM_{g,n},\mathbb{C}) \otimes (V^*)^{\otimes n}
$$
satisfying (i), (ii), and (iii) is called a {\em Cohomological Field Theory} (CohFT) {\em with unit}.
\end{definition}

The simplest example of a Cohomological Field Theory with unit is given by the
{\em trivial CohFT},
$$V=V_0 = \mathbb{Q}\, ,\ \ \eta(1,1)=1\,, \ \ \b1=1 \,, \ \
\Omega_{g,n}(1,\ldots,1) = 1\in H^0(\oM_{g,n},\mathbb{C})\, .$$
A more interesting example is given by the total Chern class 
$$c(\mathbb{E}) = 1 +\lambda_1+\ldots+ \lambda_g\, \in \, H^*(\oM_{g,n},\mathbb{C})$$
of
the rank $g$ Hodge bundle $\mathbb{E}\rightarrow \oM_{g,n}$,
$$V=V_0=\mathbb{Q}\, ,\ \ \eta(1,1)=1 \,, \ \ \b1=1\,, \ \
\Omega_{g,n}(1,\ldots,1) = c(\mathbb{E})\in H^*(\oM_{g,n},\mathbb{C})\, .$$

\begin{definition}\label{defcohft2}  For a CohFT
$\Omega= (\Omega_{g,n})_{2g-2+n>0}$, the {\em topological part}
$\omega$ of $\Omega$ is defined by 
$$\omega_{g,n} = [\Omega_{g,n}]^0 \, \in \, H^0(\oM_{g,n},\mathbb{C}) \otimes
(V^*)^{\otimes n}
\, .$$
\end{definition}

The topological (or degree 0) part $[\,\, ]^0$ of $\Omega$ is simply obtained from the 
canonical summand projection
$$[\,\, ]^0:H^*(\oM_{g,n},\mathbb{C}) \rightarrow H^0(\oM_{g,n},\mathbb{C})\, .$$
If $\Omega$ is a CohFT with unit, then $\omega$ is also a
CohFT with unit.
The topological part of the CohFT obtained from the total
Chern class of the Hodge bundle is the trivial CohFT.

A {\em Topological Field Theory} (TopFT) with unit is a CohFT 
$\theta$
with unit
of the form 
$$\Big\{\,  \theta_{g,n}\in H^0(\oM_{g,n},\mathbb{C}) \otimes (V^*)^n\, \Big\}_{2g-2+n>0}\, .$$
The topological part of a CohFT with unit  is
a TopFT with unit.

\subsection{Tautological cohomology}
The subrings of {\em tautological classes} on the moduli spaces of curves,
$$RH^*(\oM_{g,n},\mathbb{C}) \subset H^*(\oM_{g,n}, \mathbb{C})\, ,$$ have been
extensively studied --- see \cite{FPan,GrPan,Janda2,PPZ} for definitions, generators, and Pixton's conjectured
set of relations.

\begin{definition}  A CohFT
$\Omega$ {\em takes the value} $\gamma\in H^*(\oM_{g,n},\mathbb{C})$
if $\gamma$ lies in the image of
$$\Omega_{g,n}: V^{\otimes n} \rightarrow  H^*(\oM_{g,n},\mathbb{C})\, .$$
\end{definition}

As a consequence of the Givental-Teleman classification \cite{Givental,Givental2,Teleman},
{\em all} semisimple CohFTs with unit take values in the tautological cohomology
of the moduli spaces of curves.

\subsection{Constructions}
There are several constructions of CohFTs with unit --- Gromov-Witten theory,
Witten's r-spin class, and the Chern characters of the Verlinde bundle all define
CohFTs with unit, see \cite{MOPPZ,PPP,PPZ,PPZ2}. Moreover, once a CohFT is found, others can
be constructed via the action of the Givental group, see \cite{PPZ, Shadrin, Teleman}.

Our perspective here is different. We would like to construct CohFTs {\em by hand}.
The task is difficult since an infinite amount of compatibility (involving all higher
genera) is required. Our first result  concerns minimal classes.

\begin{definition}
A class $\gamma \in H^*(\oM_{g,n},\mathbb{C})$ is {\em minimal} if 
$$q^*(\gamma)=0\ \ \ \text{and} \ \ \ r^*(\gamma)=0$$
for all boundary maps to $\oM_{g,n}$ of type $q$ and $r$.
\end{definition}

In other words, $\gamma\in H^*(\oM_{g,n},\mathbb{C})$ is minimal if the
restriction of $\gamma$ to every boundary component 
of $\oM_{g,n}$ vanishes. For example,
the Poincar\'e dual of a point  in $H^{2(3g-3+n)}(\oM_{g,n},\mathbb{C})$ is
always a minimal class.

\begin{definition}
A minimal class $\gamma \in H^*(\oM_{g,n},\mathbb{C})$ satisfies the {\em parity condition} if the cohomological degree of $\gamma$ and the number $n$ of marked points have the same parity.
\end{definition}

\begin{theorem}\label{ttttt}
Let $\gamma\in H^*(\oM_{g,n},\mathbb{C})$ be a minimal class that satisfies the parity condition. Then there exists a CohFT with unit $\Omega^\gamma$
which takes the value $\gamma$.
\end{theorem}

\noindent More precisely, for every minimal class $\gamma$ satisfing the parity condition, we construct a canonical CohFT $\Omega^\gamma$ taking the value $\gamma$.

\bigskip

There exist non-tautological cohomology classes on the moduli space of curves. The simplest
is perhaps 
$$0\neq \phi \in H^{11,0}(\oM_{1,11},\mathbb{C}) \stackrel{\sim}{=} \mathbb{C}\, 
$$
defined via the discriminant modular form, see \cite[Section 2]{FPan} for
an exposition.
Since $\phi$ is a class of odd cohomological degree, 
$$\phi \notin RH^*(\oM_{1,11},\mathbb{Q})\, .$$
Since {\em no} boundary component of $\oM_{1,11}$ has nonvanishing odd cohomology,
$\phi$ is a minimal class (which also satisfies the parity condition).

\begin{corollary}
The CohFT with unit $\Omega^\phi$ takes values outside of the tautological cohomology of
the moduli spaces of curves.
\end{corollary}

The CohFT $\Omega^\phi$ is the first known example of a CohFT with unit taking non-tautological
values. 
Whether the Gromov-Witten theory of a nonsingular
projective variety $X$ can ever take non-tautological values is an interesting
question. For semisimple $X$, the Gromov-Witten CohFT must take
values in the tautological ring. Perhaps the simplest non-semisimple variety $X$
is a curve of higher genus. However, the CohFTs obtained from the Gromov-Witten
theories of target curves have been proven to take values in tautological
cohomology by Janda \cite{Jan}.

On the other hand, the Gromov-Witten theory of a nonsingular projective variety $X$ may produce classes outside of the tautological ring 
$$R^*(\oM_{g,n},\mathbb{C}) \subset A^*(\oM_{g,n},\mathbb{C})\, $$
in Chow.
Simple examples can be found in the case of 
 higher genus target curves~$X$. The
class of the moduli point{\footnote{Moduli points are
known not to be always tautological in $\oM_{g,n}$. For
example, since $\oM_{1,11}$ has a holomorphic differential
form, moduli points of $\oM_{1,11}$ are not always tautological
by \cite{Sri}.}}
$$[X,p_1,\ldots, p_n] \in A_0(\oM_{g,n}, \mathbb{C})$$
occurs as push-forward to $\oM_{g,n}$ of 
$$\prod_{i=1}^n \text{ev}_i^*(p_i) \, \cap \, \big[\oM_{g,n}(X,[X])\big]^{vir} \in A_0(\oM_{g,n}(X,[X]), \mathbb{C})\, . $$
However, the virtual class has the possibility of being
better behaved in cohomology.

\subsection{Minimal classes of even degree}
Minimal classes have played an important role in the study of
the tautological ring. For example, the tautological class,
$$\lambda_g \lambda_{g-1} \in H^{4g-2}(\MM_{g}, \com)\, ,$$
which appears in the socle evaluation of $RH^*(\cM_g)$ for $g\geq 2$,
is well-known to be minimal. 

While we expect the existence of {\em non-tautological} minimal classes of
even cohomological degree, we do not know any examples at the moment. 
Since all of the even degree cohomology in genus 1
is tautological \cite{genusone}, non-tautological minimal classes of
even degree do {\em not} exist on $\MM_{1,n}$.
D. Petersen has 
provided a proof of
the non-existence of non-tautological minimal classes of
even degree on $\MM_{2,n}$ which appears
in the Appendix. So the search for
non-tautological minimal classes of even
 degree should start in genus 3.
 
\subsection{Acknowledgments}
The construction presented here was found in Paris (at Caf\'e Nuance) 
in preparation for
the colloquium of R.P. at Jussieu in December 2017.
The results were also discussed at the conference {\em Hurwitz cycles
on the moduli of curves} at Humboldt University in Berlin in
February 2018 (organized with G. Farkas). We thank J. Schmitt for discussions
about minimal classes and D. Petersen for contributing
the Appendix on minimal  classes in low genus.

R.P. was partially supported by
 SNF grant 200021-143274,   ERC grant
AdG-320368-MCSK, SwissMAP, and the Einstein Stiftung.
D.Z. was partially supported by the grant ANR-18-CE40-0009 (ENUMGEOM).

The project has received funding from the European Research Council (ERC)
under the European Union Horizon 2020 research and innovation programme
(grant No. 786580).

\section{Construction of $\Omega^\gamma$}

Let $\gamma\in H^*(\oM_{h,m},\mathbb{C})$ be a minimal class satisfying the parity condition. The parity condition implies $(h,m) \neq (0,3)$.
Then, since $\gamma$ is minimal, the cohomological degree of $\gamma$
must be positive,
$$\gamma \in H^{>0}(\oM_{h,m},\mathbb{C})\,, \ \ \ \ (h,m) \neq (0,3)\, .$$

To construct a canonical CohFT with unit $\Omega^\gamma$ which takes the value $\gamma$, we start with the topological field theory $\omega^\gamma$ associated with the Frobenius algebra structure of $H^*(X,\C)$, where $X$ is a genus~$m$ curve{\footnote{The TopFT $\omega^\gamma$ depends
{\em only} on $m$ --- the number of markings of the moduli space
associated to $\gamma$. In Section \ref{nnnd}, we will use the notation
$\omega^m$ for $\omega^\gamma$.}}. We then modify $\omega^\gamma$ by hand by adding higher degree classes to obtain $\Omega^\gamma$. We begin with a careful description of $\omega^\gamma$.

\subsection{State space}\label{stsp}

The state space $(V,\eta,\b1)$ of both $\Omega^\gamma$ and $\omega^\gamma$ is described as follows.
\begin{enumerate}
\item[$\bullet$] Let $V$ be the $\Z_2$-graded $\C$-vector space of dimension $2m+2$ with basis given by the vectors
\begin{equation}\label{ff455}
a,b_1, b_2, \ldots, b_m, c_1,c_2, \ldots, c_m, d\,.
\end{equation}
and grading $\overline{a}=\overline{d}=0 , \; \overline{b_i} = \overline{c_i} = 1 $.
\item[$\bullet$] Let $\eta$ be the non-degenerate graded-symmetric 2-form on $V$ defined by
$$\eta(a,d)= \eta(d,a)=1\, , \;\; \eta(b_i,c_i) = -\eta(c_i,b_i)=1\, ,$$
and $\eta$ vanishes on all other pairs of basis vectors.
\item[$\bullet$] Let $\b1=a$.
\end{enumerate}
Let $\mathcal{B}$ be the set of $2m+2$ basis vectors \eqref{ff455}.
The span of  
$$a,b_1, \ldots, b_m$$
in $V$ is a maximal isotropic subspace with respect to $\eta$.
We will use the notation
$$\eta(v_1,v_2) = \langle v_1,v_2 \rangle\, .$$
The bi-vector dual to $\eta$ can be written explicitly as
\begin{equation} \label{bibi}
a \otimes d + d \otimes a - \sum_{i=1}^m (b_i \otimes c_i - c_i \otimes b_i)\, .
\end{equation}

\subsection{The algebra structure} \label{alg}

The vector space $V$ carries an algebra structure given by the following multiplication rules.
\begin{itemize}
    \item $a$ is the (left and right) unit of the algebra,
    \item $b_i \star c_i = - c_i \star b_i = d$,
    \item all pairwise products of basis elements vanish, except in the two cases above.
\end{itemize}


Because the above algebra is isomorphic to the cohomology algebra of a genus~$m$ surface, there is a natural $\Z$-grading which lifts the $\Z_2$-grading:
$$
\overline{a} = 0\, , \quad \overline{b_i} = \overline{c_i} = 1\, , \quad \overline{d}=2\, .
$$

\begin{remark}
The algebra structure on~$V$ is {\em not} semisimple. Indeed, all elements of the basis $\mathcal{B}$ except $a$ are nilpotents.
\end{remark}

\subsection{The values of $\omega^\gamma$}

Since $\omega_{g,n}^\gamma$ takes values in
$H^0(\cM_{g,n},\mathbb{C})$ which is canonically
$\mathbb{C}$, we view the values 
of $\omega_{g,n}^\gamma$
as complex numbers.

\begin{proposition} \label{p2p2}
 The TopFT $\omega^\gamma$
 has the following evaluation on basis elements of $\mathcal{B}$.
\begin{enumerate}
\item[$\bullet$] In genus $g=0$,
$\omega^\gamma_{0,n}(d, \underbrace{a, \ldots,a}_{n-1}) = 1$ and 
$$\omega^\gamma_{0,n}(b_i, c_i \underbrace{a, \ldots,a}_{n-2}) = 1\, \ \ 
{\text{for}}\ \ 1\leq i \leq m\, .$$ The evaluations obtained from these by permuting the entries are equal to $1$ or $-1$ as determined by the grading. All other evaluations vanish.
\item[$\bullet$] In genus $g=1$,
$\omega^\gamma_{1,n}(a,\ldots,a) = 2-2m$
and all other evaluations vanish.
\item[$\bullet$] In genus $g\geq 2$, all evaluations vanish,
 $\omega^\gamma_{g,n}= 0$.
\end{enumerate}
\end{proposition}

\paragraph{Proof.}
By the axioms of a TopFT, the genus~0 evaluation $\omega^\gamma_{0,n} (v_1, \dots, v_n)$ is equal to $\eta( v_1 \star \dots  \star v_n, \b1)$, which is the coefficient of $d$ in the product $ v_1 \star \dots  \star v_n$. Using the $\Z$-grading of the algebra, the coefficient vanishes unless the sum of the 
gradings of $v_1, \dots, v_n$ equals~2, whence the result.

For the higher genus cases, we evaluate the $0$-cohomology class $\omega^\gamma_{g,n}(v_1, \dots, v_n)$ on a point of the moduli space that corresponds to a genus~$g$ curve with $g$ nonseparating nodes, in other words, a rational curve with $g$ pairs of identified points. Using the genus~0 case, only one node is possible without vanishing, therefore $\omega^\gamma$ vanishes for $g \geq 2$. The factor $2-2m$ (the Euler characteristic of a genus~$m$ curve) in the $g=1$ case occurs via CohFT axiom (ii) and the definition of the symmetric form~$\eta$ and its inverse bi-vector.   
\qed

\subsection{Full CohFT $\Omega^\gamma$}
Recall $\gamma\in H^*(\oM_{h,m},\mathbb{C})$ is a minimal class satisfying
$$\gamma \in H^{>0}(\oM_{h,m},\mathbb{C})\,, \ \ \ \ (h,m) \neq (0,3)\, .$$ 
We also assume that $\gamma$ satisfies the parity condition:
deg $\gamma$ has the same parity as the number of points $m$.

\begin{definition} \label{kkkk}
The values of the full CohFT $(\Omega^\gamma_{g,n})_{2g-2+n>0}$ on basis vectors from $\mathcal{B}$ are defined by the following rules.
\begin{enumerate}
\item[$\bullet$] Let $g=h$ and $n\geq m$. Let $p : \oM_{h,n} \to \oM_{h,m}$ be the forgetful map.
Then we have
\begin{equation} \label{Eq1}
\Omega^\gamma_{g,n}(b_1, b_2, \dots b_m, \underbrace{a,\ldots,a}_{n-m}) 
= p^*\gamma\, .
\end{equation}
\item[$\bullet$] The evaluations obtained from the above by permutations of the entries are equal to $p^*\gamma$ or $-p^*\gamma$ according to the grading.
\item[$\bullet$] In all other cases, namely, either $g \not= h$, or $g=h$, but $n<m$, or $g=h$ and $n \geq m$, but the basis vectors on which $\Omega^\gamma$ is evaluated are not obtained by a permutation of~\eqref{Eq1}, we have
$$\Omega^\gamma_{g,n}(v_1,\ldots,v_n)= \omega^\gamma_{g,n}(v_1,\ldots,v_n)\, .$$
\end{enumerate}
\end{definition}

The parity condition $\text{deg}\, \gamma=m \ \text{mod}\ 2$ implies that the
tensor $\Omega^\gamma_{g,n}$ is even.
Certainly $\Omega^\gamma$ takes the value $\gamma$ since, by definition,
$$\Omega^\gamma_{h,m}(b_1,\ldots,b_m)= \gamma \in H^*(\oM_{h,m},\mathbb{C})\, .$$
To complete the proof of Theorem \ref{ttttt}, we must check 
that $\Omega^\gamma$
satisfies all of the required axioms for a CohFT with unit.
CohFT axiom (i), 
invariance under the symmetric group, and CohFT axiom (iii),
compatibility with the forgetful map $p$,
both follow immediately from the construction.

\subsection{CohFT axiom (ii) for $\Omega^\gamma$}
Our first remark is that \eqref{Eq1} of Definition~\ref{kkkk} can be rewritten equivalently as
$$
\Omega^\gamma_{g,n}(b_1,\ldots,b_m, \underbrace{a, \dots, a}_{n-m}) = \omega^\gamma_{g,n}(b_1,\ldots,b_m, \underbrace{a, \dots, a}_{n-m}) + p^*\gamma.
$$
Indeed, $\omega^\gamma_{g,n}(b_1,\ldots,b_m, \underbrace{a, \dots, a}_{n-m})=0$ for all $g$ and~$n$ as shown in Proposition~\ref{p2p2}. (If $g \geq 2$ then all correlators vanish. If $g=0$ we need either a $d$-insertion or a $b$- and a $c$-insertion to get a nonvanishing correlator, which we don't have here. Finally, if  $g=1$  a nonvanishing correlator has  $a$-insertions only. This could be the case here if $m=0$. But since $\gamma$ is a cohomology class in $\oM_{h,m}$, the case $h=1, m=0$ is ruled out.)

In order to check CohFT axiom (ii) for 
$\Omega^\gamma$, we must prove the compatibility of $\Omega^\gamma$ under
pull-back for the boundary morphisms of type $q$ and $r$.
Since CohFT axiom (ii) holds for $\omega^\gamma$, we must only
study the effect of the correction term $p^*(\gamma)$.

Consider first the boundary morphism 
$$
q : \oM_{g-1, n+2} \to \oM_{g,n}\, .
$$
There are two cases where the correction term has the possibility of
appearing:
\begin{enumerate}
\item[$\bullet$] If $g-1=h$, the correction term could appear in the left-hand side. However CohFT axiom~(ii) with the bi-vector~\eqref{bibi}
includes basis elements of type $c$ or~$d$ at one of the last two markings, while the correction term only appears when all basis elements are of types $a$ and~$b$. So the correction never appears and the compatibility holds. 
\item[$\bullet$] If $g=h$, $n\geq m$, the correction term appears in the right-hand side whenever the markings $(v_1,\ldots,v_n)$ are a permutation of
$$(b_1,b_2,\ldots,b_m, \underbrace{a,\ldots,a}_{n-m})\,.
$$
Then $q^*(\Omega^\gamma_{h,n}(v_1,\ldots,v_n))$ vanishes since $\gamma$ is minimal. Thus the minimality of $\gamma$ ensures that the compatibility still holds.
\end{enumerate}

Consider next the boundary morphism
$$
r: \oM_{g_1, n_1+1} \times \oM_{g_2, n_2+1} \to \oM_{g,n}\, , 
$$
where $n=n_1+n_2$ and $g=g_1+g_2$. 

On the right-hand side, the correction term appears in
$$
\Omega^\gamma_{h,n}(b_1, \dots, b_m, \underbrace{a, \dots, a}_{n-m})$$
and terms obtained from this by permutations of entries. We
 enumerate the cases where the correction term appears on 
the left-hand side as follows. 
Using Proposition~\ref{p2p2}, there are exactly four cases, 
up to permuting the $n_1$ entry vectors of the first factor and the $n_2$ entry vectors in the second factor. 
In the formulas below, a vector with a hat denotes a skipped entry, 
while an underlined vector originates from the sum $\sum \eta^{jk} e_j \otimes e_k$ :
\begin{itemize}
    \item $\Omega^\gamma_{h,n_1+1}(b_1, \dots, b_m, \underbrace{a, \dots, a}_{n_1-m}, \underline{a}) \times
    \Omega^\gamma_{0,n_2+1}(\underline{d}, \underbrace{a, \dots, a}_{n_2})$,
     \item $\Omega^\gamma_{h,n_1+1}(b_1, \dots, \widehat{b}_i, \dots, b_m, \underbrace{a, \dots, a}_{n_1-m+1}, \underline{b_i}) \times
    \Omega^\gamma_{0,n_2+1}(\underline{c_i}, b_i, \underbrace{a, \dots, a}_{n_2-1})$,
    \item $\Omega^\gamma_{0,n_1+1}(\underbrace{a, \dots, a}_{n_1},\underline{d}) \times \Omega^\gamma_{h,n_2+1}(\underline{a},b_1, \dots, b_m, \underbrace{a, \dots, a}_{n_2-m}),$
     \item $\Omega^\gamma_{0,n_1+1}(b_i, \underbrace{a, \dots, a}_{n_1-1},\underline{c_i}, ) \times\Omega^\gamma_{h,n_2+1}(\underline{b_i}, b_1, \dots, \widehat{b}_i, \dots, b_m, \underbrace{a, \dots, a}_{n_2-m+1})$.
\end{itemize}
To show that these are the only cases, we use two simple remarks: 
\begin{enumerate}
\item[(i)]Since each term in the bi-vector $\eta^{-1}$ always contains a basis vector of type $c$ or $d$ as a factor, the correction term can appear only in one of the two factors.
\item[(ii)]  The factor without correction contains an entry of $c$ or $d$ type and thus can only be of genus~0 so as not to vanish.
\end{enumerate}

We can now prove the equality 
\begin{multline*}
r^*(\Omega_{g,n}(v_1,\ldots,v_n)) = \\
\sum_{j,k} \eta^{jk}\, \Omega_{g_1,n_1+1}(v_1,\ldots,v_{n_1}, e_j) \otimes
\Omega_{g_2,n_2+1}(e_k, v_{n_1+1}, \ldots, v_n)\, .
\end{multline*}
First, from the analysis above, we see that both sides vanish unless $g=h$ and $v_1, \dots, v_n$ is a permutation of $b_1, \dots, b_m, \underbrace{a, \dots, a}_{n-m}$. 

Now assume  $g=h$ and $v_1, \dots, v_n$ is a permutation of $b_1, \dots, b_m, \underbrace{a, \dots, a}_{n-m}$. Let $p : \oM_{h,n} \to \oM_{h,m}$ be the  map which
 forgets the marked points carrying the basis vector~$a$. Consider the image of the map $p \circ r$. We have 
$$
r^* \Omega_{g,n}(v_1, \dots, v_n) = \pm (p \circ r)^* \gamma.
$$
By the minimality of $\gamma$, this class vanishes whenever the image of $p \circ r$ is a boundary stratum. In order for $p \circ r$ to be onto, $p$ has to contract one of the two irreducible components of the curve. A component is contracted if it has genus~0 and all of its marked points, except perhaps one, carry the basis vector~$a$. This leaves us with exactly the same four cases as in the enumeration above. In the first two cases, both 
$$
r^* \Omega_{g,n}(v_1, \dots, v_n)
$$
and
$$
\sum _{j,k} \eta^{jk} \Omega_{g_1,n_1+1}(v_1, \dots, v_{n_1}, e_j)
 \otimes
 \Omega_{g_2,n_2+1}(e_k, v_{n_1+1}, \dots v_n)
$$
are equal to
$$
\pm (p \circ r)^* \gamma \otimes 1\, ,
$$
the sign determined by the order of the vectors, which is the same on both sides. Similarly, in the last two cases, 
both 
$$
r^* \Omega_{g,n}(v_1, \dots, v_n)
$$
and 
$$
\sum _{j,k} \eta^{jk} \Omega_{g_1,n_1+1}(v_1, \dots, v_{n_1}, e_j)
 \otimes
 \Omega_{g_2,n_2+1}(e_k, v_{n_1+1}, \dots v_n)
$$
are equal to
$$
\pm 1 \otimes (p \circ r)^* \gamma\, .
$$
Again, the sign is determined by the order of the vectors, which is the same on both sides. The proof of Theorem \ref{ttttt} is complete. \qed

\begin{remark}
{\em  The strategy can be
summarized as follows. 
Given a minimal class
$$\gamma\in H^{>0}(\oM_{h,m},\mathbb{C})\, ,$$
we start with a Topological Field theory
$\omega^\gamma$ which depends only on the number of markings $m$.
We then define $\Omega^\gamma$ by adding a correction term to
$\omega^\gamma$ in the particular cases 
\begin{equation} \label{dxxd}
\Omega^\gamma_{h,n}(b_1,\ldots,b_m,\underbrace{a,\ldots, a}_{n-m}) = \omega^\gamma_{h,n}(b_1,\ldots,b_m,\underbrace{a,\ldots, a}_{n-m})
\ +\ p_{n,m}^*(\gamma)\, ,
\end{equation}
up to permutations of the entries. To verify the CohFT axioms for $\Omega^\gamma$, we use the
minimality of $\gamma$ for lower moduli spaces and the
exact solution of the TopFT $\omega^\gamma$ of Proposition \ref{p2p2}
for higher moduli spaces. It is important that the correction~\eqref{dxxd} occurs only for insertions in the maximal isotropic
subspace spanned by 
$a,b_1,\ldots, b_m$.} 
\end{remark}

\begin{remark}
{\em If $\gamma$ is a minimal class of even degree, it is possible to remove the grading from the construction. The basis of the vector space~$V$ remains the same, but the quadratic form $\eta$ becomes symmetric and, in general, all the signs related to the grading disappear. In this case, it is not necessary to require the number of marked points $m$ to have the same parity as $\gamma$.} 
\end{remark}

\section{Deformations}
\subsection{General theory}\label{gent}
Let $\Omega$ be a CohFT with unit and state space $(V,\eta,\b1)$.

\begin{definition} A system of tensors
$$\Lambda=\Big(\, \Lambda_{g,n}:V^{\otimes n} \rightarrow H^*(\oM_{g,n},\mathbb{C})\,  \Big)_{2g-2+n>0}$$
defines a {\em first order deformation  of $\Omega$}
if $$\Big(\, \Omega_{g,n}+ \epsilon \Lambda_{g,n}\, \Big)_{2g-2+n>0}$$
satisfies the axioms of a  CohFT with unit modulo $\epsilon^2=0$.
\end{definition}

If the system of tensors $\Lambda=(\Lambda_{g,n})_{2g-2+n>0}$ satisfies the further condition
$$\Lambda_{0,3} =0\,, $$
then the first order deformation $\Lambda$ of $\Omega$  preserves the 
TopFT structure.
The proof of Theorem~\ref{ttttt} shows that every minimal class
$$\gamma \in H^{>0}(\oM_{h,m},\mathbb{C})$$
whose degree has the same parity as~$m$
yields a first order deformation of $\omega^\gamma$ which preserves the TopFT structure.

We can write the CohFT axiom conditions for the deformation
 $$\Big(\, \Omega_{g,n}+ \epsilon \Lambda_{g,n}\, \Big)_{2g-2+n>0}$$
more
explicitly:
\begin{enumerate}
\item[(i)] Each tensor $\Lambda_{g,n}$ is $S_n$-invariant(in
the $\mathbb{Z}_2$-graded sense) for
 the natural action of the symmetric group $S_n$.
\item[(iiq)] 
The tensor $q^*(\Lambda_{g,n}) \in H^*(\oM_{g-1,n+2},\mathbb{C}) \otimes 
(V^*)^{\otimes n}$,
obtained via pull-back by the boundary morphism 
$$
q : \oM_{g-1, n+2} \to \oM_{g,n}\ ,
$$
is required to equal the contraction
of $\Lambda_{g-1,n+2}$ by the bi-vector 
$$\sum_{j,k} \eta^{jk} e_j \otimes e_k\, .$$
\item[(iir)]
The tensor $r^*(\Lambda_{g,n})$, obtained via pull-back 
by the boundary morphism
$$
r: \oM_{g_1, n_1+1} \times \oM_{g_2, n_2+1} \to \oM_{g,n}\, , 
$$
is required
to equal 
\begin{eqnarray*}
& &\sum_{j,k} \eta^{jk}\, \Omega_{g_1,n_1+1}(v_1,\ldots,v_{n_1}, e_j) \otimes
\Lambda_{g_2,n_2+1}(e_k,v_{n_1+1}, \ldots, v_n) \\
& + & \sum_{j,k} \eta^{jk}\, \Lambda_{g_1,n_1+1}(v_1,\ldots,v_{n_1}, e_j) \otimes
\Omega_{g_2,n_2+1}(e_k,v_{n_1+1}, \ldots, v_n)\,.
 \end{eqnarray*}

\item[(iii)] The tensor $p^*(\Lambda_{g,n})$, obtained via pull-back by the forgetful map
$$p:\oM_{g,n+1} \rightarrow \oM_{g,n}\, ,$$
is required to satisfy 
$$
\Lambda_{g,n+1}(v_1, \ldots,v_n , \b1) = p^*\Lambda_{g,n} (v_1, \ldots,  v_n)\ 
$$
for all $v_i \in V$.
\end{enumerate}

The simplest method of constructing deformations of $\Omega$
is via Givental's $R$-matrix action. Other deformations are,
in general, hard to find.

\subsection{Isotropic deformations of the CohFT $\omega^m$}\label{nnnd}
Consider the topological field theory $\omega^m$ defined in Sections~\ref{stsp} and~\ref{alg}. The state space $V$ has basis
\begin{equation*}
\mathcal{B} = (a,b_1, b_2, \ldots, b_m, c_1,c_2, \ldots, c_m, d)\,,
\end{equation*}
and carries the non-degenerate graded-symmetric
2-form on $\eta$ defined by
$$\eta(a,d)=\eta(d,a)=1, \quad \eta(b_i, c_i) =- \eta(c_i, b_i) =  1$$
with values on all other pairs of basis vectors defined to vanish.

\begin{definition} 
A first order deformation of $\omega^m$ 
defined by $\Lambda$ is {\em isotropic} if 
$$\Lambda_{g,n}(v_1,\ldots,v_{n-1},c_i)= \Lambda_{g,n}(v_1,\ldots,v_{n-1},d)= 0$$
for all $g$, $n$, $v_j\in \mathcal{B}$, and $1\leq i \leq n$.
\end{definition}

\begin{theorem} \label{ww22}  Let $\Lambda$ define
an isotropic first order deformation of $\omega^m$ 
 which preserves the TopFT structure. Then,
$$\Lambda_{g,n}(v_1,\ldots,v_n)\in  H^*(\oM_{g,n},\mathbb{C})\ \ \ \text{for} \ \ \ v_j \in \{b_1,\ldots,b_m\}$$
is always a minimal class.
\end{theorem}

\paragraph{Proof.}
Let $(v_1,\ldots,v_n)$ be an $n$-tuple of vectors satisfying
$v_j \in \{b_1,\ldots,b_m\}$ for all $j$.
In order to prove that the class
$$\Lambda_{g,n}(v_1,\ldots,v_n)\in  H^*(\oM_{g,n},\mathbb{C})$$
is minimal, we must study the pull-backs under the boundary
morphisms of type $q$ and $r$.
\begin{enumerate}
\item[$\bullet$]
The vanishing of $q^*(\Lambda_{g,n}(v_1,\ldots,v_n))$ under
$$
q : \oM_{g-1, n+2} \to \oM_{g,n}\, 
$$
is immediate from axiom (iiq) of Section \ref{gent}:
each term of the bi-vector~\eqref{bibi}
always includes a basis element of type~$c$ or~$d$.
\item[$\bullet$] The vanishing of
$r^*(\Lambda_{g,n}(v_1,\ldots,v_n))$ under
$$
r: \oM_{g_1, n_1+1} \times \oM_{g_2, n_2+1} \to \oM_{g,n}\, , 
$$
where $n=n_1+n_2$ and $g=g_1+g_2$, follows
directly from axiom (iir), Proposition \ref{p2p2} for $\omega^m$, and
the isotropic property of $\Lambda$.
\end{enumerate}
Since $\Lambda_{g,n}(v_1,\ldots,v_n)$ vanishes under all boundary
restrictions, the class is minimal.
\qed
\vspace{8pt}

Theorems~\ref{ttttt} and~\ref{ww22} show minimal
classes $\gamma$ satisfying $\deg \gamma = m \bmod 2$ and the isotropic deformations of $\omega^m$ are
essentially equivalent notions. 
Deformations of $\omega^m$ via Givental's $R$-matrix action will
in general not be isotropic. Also, the $R$-matrix action
 will produce only deformations
taking values in the tautological cohomology --- so the deformations
obtained from minimal classes $$\gamma \in H^{>0}(\oM_{h,m},\mathbb{C})$$
will not always be obtained via the $R$-action. 

\pagebreak

\begin{center} {\Large{{\noindent{\bf {\sf Appendix:
{Minimal cohomology classes on $\MM_{g,n}$\\ in low genus}}}}}}
\end{center}

\begin{center}{\large {by D. Petersen}} \end{center}

\setcounter{section}{0}
\noindent{\large{\bf {A1. Minimal cohomology}}}

\vspace{8pt}
In this appendix, we make some remarks regarding minimal cohomology classes on $\MM_{g,n}$. Cohomology will always be taken with $\Q$-coefficients. We begin by recalling a useful result from Deligne's mixed Hodge theory.

\begin{prop}\label{deligne-lemma}
	Let $X$ be a smooth projective variety, $Z \subset X$ a closed subvariety{\footnote{$Z$ need not be irreducible.}} of pure codimension $c$ with complement $U$, $\widetilde Z \to Z$ any resolution of singularities. There are exact sequences
	$$ H^{k-2c}(\widetilde Z)(-c) \to H^k(X) \to \mathrm{Gr}^W_k H^k(U) \to 0 $$
	and 
	$$ 0 \to \mathrm{Gr}^W_k H^k_c(U) \to H^k(X) \to H^k(\widetilde Z).$$
\end{prop}

\paragraph{Proof.}
The first exact sequence follows by combining \cite[Corollaire 3.2.17]{hodge2} and \cite[Corollaire 8.2.8]{hodge3}. The second one is just the Poincar\'e dual of the first one.
\qed \vspace{8pt}

We apply Proposition \ref{deligne-lemma} to the case where $X = \MM_{g,n}$, $U = \cM_{g,n}$, and $\widetilde Z$ is the normalization of the boundary, so that each component{\footnote{For some components, $\widetilde Z$ is $2$-fold
cover of the normalization, but Proposition 1 still holds.}}
 of $\widetilde Z$ is either of the form $\MM_{g',n'+1} \times \MM_{g-g',n-n'+1}$ or $\MM_{g-1,n+2}$.

\begin{defn}Let $H^k_{\min}(\MM_{g,n}) \subseteq H^k(\MM_{g,n})$ denote the subspace of minimal cohomology classes.
\end{defn}

\begin{prop} 
		For all $k$, $g$, and $n$, there is an isomorphism 
		$$ \mathrm{Gr}^W_k H^k_c(\cM_{g,n}) \cong H^k_{\min}(\MM_{g,n})\, .$$
\end{prop}

\paragraph{Proof.}
Immediate from the second exact sequence of Proposition \ref{deligne-lemma}.
\qed \vspace{8pt}

\begin{prop} \label{perfectpairing}
	There is a perfect pairing between 
$$\mathrm{Gr}^W_k H^k(\cM_{g,n}) \ \ \text{and} \ \ H^{2(3g-3+n)-k}_{\min}(\MM_{g,n})\, .$$
\end{prop} 

\paragraph{Proof.}
There is a perfect pairing between $H^k(\cM_{g,n})$ and $H^{2(3g-3+n)-k}_c(\cM_{g,n})$ by Poincar\'e duality, which induces a perfect pairing between their respective associated graded for the weight filtration. With the previous proposition, the result follows.
\qed \vspace{8pt}

In fact the perfect pairing of Proposition \ref{perfectpairing} can be described more explicitly. Take any class $\alpha \in \mathrm{Gr}^W_k H^k(\cM_{g,n})$, and lift it (non-canonically) to a class $\overline \alpha \in H^k(\MM_{g,n})$. If $\beta \in H^{2(3g-3+n)-k}_{\min}(\MM_{g,n})$ is a class of complementary degree, then the cup product
$$ \overline{\alpha} \cup \beta \in H^{2(3g-3+n)}(\MM_{g,n}) \cong \Q$$
is in fact well defined: if $\overline \alpha '$ is a different choice of lift, then the difference $\overline \alpha - \overline \alpha'$ is pushed forward from the boundary, so its integral against $\beta$ is zero since $\beta$ vanishes on the boundary. The reader familiar with the $\lambda_g \lambda_{g-1}$-pairing on the tautological ring of $\cM_g$ will find this construction familiar; indeed, what makes the $\lambda_g\lambda_{g-1}$-pairing work is precisely that $\lambda_g\lambda_{g-1}$ is a minimal class on $\MM_g$.

\vspace{12pt}
\noindent{\large{\bf{A2. Minimal classes in genus zero}}}

\begin{prop}\label{genuszero}
	The point class in $H^{2(n-3)}(\MM_{0,n})$ is the only minimal class in $H^*(\MM_{0,n})$. 

\end{prop}

\paragraph{Proof.}
There are very many ways to see this; here is one. By Proposition \ref{perfectpairing}, the claim is equivalent to  $\mathrm{Gr}^W_k H^k(\cM_{0,n}) = 0$ for $k>0$. Consider the compactification $\cM_{0,n} \subset \C^{n-3} \subset \mathbf 
{\mathbb{P}}^{n-3}$. By the first exact sequence of Proposition \ref{deligne-lemma}, there is a surjection $H^k({\mathbb{P}}^{n-3}) \to \mathrm{Gr}^W_k H^k(\cM_{0,n})$. This map factors through $H^k(\C^{n-3})$, which vanishes for $k>0$.
\qed \vspace{8pt}

\vspace{4pt}
\noindent{\large{\bf A3. Minimal classes in genus one}}

\begin{prop}\label{genusone}
	The point class in $H^{2n}(\MM_{1,n})$ is the only minimal class of even degree in $H^*(\MM_{1,n})$. 
\end{prop}

\paragraph{Proof.}
	By Proposition \ref{perfectpairing}, this is equivalent to the claim that $\mathrm{Gr}^W_{2k}H^{2k}(\cM_{1,n})$ is nontrivial only for $k=0$, which is exactly \cite[Theorem 1.1]{genusone}.
\qed \vspace{8pt}

As pointed out in the body of this paper, the existence of odd cohomology in $H^*(\MM_{1,n})$ for $n \geq 11$ implies that there exist plenty of odd minimal classes in genus one. For any positive integer $k$ there is a rational Hodge structure $S[k]$ which is ``attached'' to the space $S_k$ of cusp forms for $\mathrm{SL}(2,\Z)$ of weight $k$. Over $\C$, the Hodge structure becomes canonically the direct sum of the spaces of holomorphic and antiholomorphic cusp forms:
$$ S[k] \otimes_{\Q} \C \cong S_k \oplus \overline{S_k}.$$
The Hodge structure $S[k]$ is pure of weight $k-1$, and the Hodge numbers of the two summands above are $(k-1,0)$ and $(0,k-1)$. Arguments similar to those of \cite[Section 2 and 3]{genusone} allow one to prove the following more refined statement.

\begin{prop}\label{cuspforms}There is an isomorphism 
	$$\mathrm{Gr}^W_k H^k(\cM_{1,n}) \cong S[k+1] \otimes \mathrm{Ind}_{S_k \times S_{n-k}}^{S_n} \mathrm{sgn}_k \otimes \mathbf 1_{n-k}$$
	for all $n \geq k$. When $n < k$, $\mathrm{Gr}^W_k H^k(\cM_{1,n})=0$. 
\end{prop}

\noindent We omit the proof. 

\vspace{12pt}
\noindent{\large{\bf{A4. Minimal classes from tautological classes}}}
\vspace{8pt}

Before moving on to genus two, let us point out an obvious source of classes in $\mathrm{Gr}^W_k H^k(\cM_{g,n})$: every tautological class on $\cM_{g,n}$ is of pure weight, being the class of an algebraic cycle. Thus every tautological class in $RH^k(\cM_{g,n})$ must ``give rise to'' a minimal class in degree $2(3g-3+n) - k$, where we put quotation marks to emphasize that the vector spaces $\mathrm{Gr}^W_k H^k(\cM_{g,n})$ and $H^{2(3g-3+n)-k}_{\min}(\MM_{g,n})$ are {dual} to each other and a vector in one of them does not correspond canonically to any vector in the other one.

By a result of Buryak-Shadrin-Zvonkine \cite{toptautologicalgroup}, the tautological ring of 
$\cM_{g,n}$ vanishes above cohomological degree $2(g-1)$ for $n > 0$, 
and $RH^{2(g-1)}(\cM_{g,n})$ is spanned by the classes $\psi_1^{g-1}, \psi_2^{g-1}, \ldots, \psi_n^{g-1}$. These classes should correspond to certain minimal classes in the cohomology of $\MM_{g,n}$, and this is in fact explained in \cite{toptautologicalgroup}: 
let $$\alpha_s = \lambda_g \lambda_{g-1} \psi_1 \psi_2 \dots \widehat{\psi_s} \cdots \psi_n\,$$
 where the hat means an omitted factor. In \cite[Section 2]{toptautologicalgroup} the authors show that the classes $\alpha_s$ are minimal, and that the pairing between the $n$ classes $\alpha_1,\ldots,\alpha_n$ and $\psi_1^{g-1},\ldots,\psi_n^{g-1}$ is perfect. The authors of  \cite{toptautologicalgroup} also proposed a generalization of the Faber conjecture: that the ring $RH^*(\cM_{g,n})$ is level of type $n$, i.e.\ that a class in the tautological ring vanishes  if and only if its product with any class of complementary degree vanishes. We expect this to be false. However, if we \emph{assume} this statement for the moment, then every nonzero $\alpha \in RH^k(\cM_{g,n})$ pairs nontrivially with a class in $RH^{2(3g-3+n)-k}_{\min}(\MM_{g,n})$, and the collection of minimal classes on $\MM_{g,n}$ that ```corresponds to'' $RH^*(\cM_{g,n})$ is just the ideal in $RH^*(\MM_{g,n})$ generated by the $n$ minimal classes $\alpha_1,\ldots,\alpha_s$.  

\begin{rem} {\em In genus zero and one, every cohomology class (resp.\ every even degree cohomology class) on $\MM_{g,n}$ is tautological. Moreover, the tautological ring of $\cM_{g,n}$ is trivial in these cases. This is  another way of seeing Propositions \ref{genuszero} and \ref{genusone}.}
\end{rem}

\pagebreak
\vspace{8pt}
\noindent{\large{\bf{A5. More interesting examples in genus two}}}
\vspace{8pt}

When $g=2$, the tautological ring $RH^*(\cM_{2,n})$ has rank $n$ in cohomological degree $2$ and vanishes above this, and, as explained in the previous subsection, one finds a subspace of $n$ minimal classes in $H^{4+2n}_{\min}(\MM_{2,n})$. However, in genus two it is not true that every class in $\mathrm{Gr}^W_{2k} H^{2k}(\cM_{2,n})$ is tautological. The first time one sees a pure even degree class on $\cM_{2,n}$ is when $n=20$, in which case one finds the class constructed by Graber--Pandharipande \cite{GrPan}: the class in $H^{22}(\MM_{2,20})$ of bielliptic curves whose markings are switched pairwise by the bielliptic involution is nontautological, and its restriction to $H^{22}(\cM_{2,20})$ is nontrivial. In fact $\mathrm{Gr}^W_{22}H^{22}(\cM_{2,20})$ is spanned by the bielliptic class and its conjugates under the $\mathbb S_{20}$-action, and equals a copy of the representation $[2,2,\ldots,2]$ of $\mathbb S_{20}$. By Proposition \ref{perfectpairing}, this shows that $H^{24}_{\min}(\MM_{2,20})$ is also spanned by a copy of the representation $[2,2,\ldots,2]$. But it is also known that the Graber--Pandharipande class (and its $\mathbb S_{20}$-conjugates) are the only nontautological even degree classes on $\MM_{2,20}$, in the sense that the span of these classes and the tautological ring is the whole even cohomology. In particular, the classes in $H^{24}_{\min}(\MM_{2,20})$ must be tautological. The assertions of this paragraph are proven in \cite{m28ct,petersentommasi}.

By what we have said so far, there must exist a tautological minimal class on $\MM_{2,20}$, very different from the ``obvious'' examples given by the classes $\alpha_s$ from the previous subsection. 

\begin{problem} Find a geometric construction{\footnote{The problem is formulated
in Chow.}} of a minimal class in $R^{12}(\MM_{2,20})$.
\end{problem}

\noindent As we will see shortly, the class must in fact be the pushforward of a class in $R^{11}(\MM_{1,22})$.

Hain and Looijenga \cite[Conjecture 5.4]{hainlooijenga} at one point conjectured that the ideal of minimal classes in the tautological ring of $\MM_g$ is principal, generated by $\lambda_g\lambda_{g-1}$; this was part of a proposed generalization of Faber's conjecture on $\cM_g$. The obvious generalization to incorporate marked points would be that the ideal of minimal classes in 
$R^*(\MM_{g,n})$ is generated by $\alpha_1,\ldots,\alpha_n$. From what we have said here, such a conjecture is false (and fails ``for the first time'' on $\MM_{2,20}$).

For $n > 20$, one finds a larger and larger number of nontautological classes, and one could hope to find a nontautological even minimal class in genus two. Unfortunately, this is not possible:
\begin{prop} Every class in $H^{2k}_{\min}(\MM_{2,n})$ is pushed forward 
along the boundary map
from $H^{2k-2}(\MM_{1,n+2})$.
In particular, every even minimal cohomology class in genus two is tautological.
\end{prop}

\paragraph{Proof.}
From Proposition \ref{deligne-lemma} we get the short exact sequence
	$$ H^{k-2}(\MM_{1,n+2}) \to H^k(\MM_{2,n}) \to H^k(\cM_{2,n}^{ct})\, . $$
	This shows that classes not pushed forward from $\MM_{1,n+2}$ map injectively into $H^k(\cM_{2,n}^{ct})$. Since the map from $\mathrm{Gr}^W_k H^k_c(\cM_{2,n}) = H^k_{\min}(\MM_{2,n})$ factors through $H^k_c(\cM_{2,n}^{ct})$, it is enough to argue that a minimal class of even degree which is in the image of $H^k_c(\cM_{2,n}^{ct}) \to H^k(\cM_{2,n}^{ct})$ is actually zero. 
	
	Now we apply the results of \cite{m28ct}. From the study of the decomposition theorem in that paper we see that if $f \colon M_{2,n}^{ct} \to M_2^{ct}$ is the forgetful map then there are isomorphisms 
	$$ H^k_c(\cM_{2,n}^{ct}) \cong \bigoplus_{p+q=k} H^p_c(\cM_2^{ct},R^qf_\ast 
\Q), \quad H^k(\cM_{2,n}^{ct}) \cong \bigoplus_{p+q=k} H^p(\cM_2^{ct},R^qf_\ast \Q),$$
	compatible with the map $H^k_c(\cM_{2,n}^{ct}) \to H^k(\cM_{2,n}^{ct})$, and that the sheaves $R^qf_\ast \Q$ decompose as direct sums of local systems on 
$\cM_2^{ct}$ and on $\mathrm{Sym}^2 \cM_{1,1}$ associated to representations of the symplectic group. 
	
	The cohomology of local systems on $\cM_2^{ct} = A_2$ is known from \cite{localsystemsA2}, where it is in particular shown that the map $H^k_c(A_2,V_\lambda) \to H^k(A_2,V_\lambda)$ can only be nontrivial in the middle degree $k=3$. In particular a minimal class on $\cM_{2,n}^{ct}$ which lies in a summand corresponding to such a local system must be of odd degree, as the sheaves $R^qf_\ast \Q$ have vanishing cohomology for odd $q$ because of the hyperelliptic involution. On the other hand the cohomologies of local systems on $\mathrm{Sym}^2 \cM_{1,1}$ will never give rise to nontrivial minimal classes; those summands of $H^\ast(\cM_{2,n}^{ct})$ restrict isomorphically to corresponding summands in the cohomology of the preimage of $\mathrm{Sym}^2 \cM_{1,1}$, i.e.\ to the cohomology of $\cM_{2,n}^{ct} \setminus \cM_{2,n}^{rt}$. \qed \vspace{8pt}

\begin{rem}
{\em	An important ingredient in the previous proof is that we understand completely the cohomology of local systems in genus two. In genus three there is only partial conjectural information based on point counts \cite{bfg11}. Assuming conjectural formulas for the cohomology of local systems in genus three, it seems plausible that one can obtain a similar result also when $g=3$: any minimal even class on $\MM_{3,n}$ is pushed forward from $\MM_{2,n+2}$. 
Since there are even non-tautological classes in genus two, this does not rule
out the existence of minimal even non-tautological classes in genus three.
The same conjectural formulas suggest that the first case where one finds nontautological classes is $\MM_{3,18}$.}
\end{rem}

\begin{rem}
{\em Just as in genus one, there are lots of odd minimal classes on $\MM_{2,n}$, which can be described in terms of automorphic forms. Although there is no statement as simple as Proposition \ref{cuspforms}, it turns out that every vector-valued cusp eigenform for $\mathrm{Sp}(4,\Z)$ of weight $\geq 3$ gives rise to cohomology classes{\footnote{For an introduction to these
constructions, see \cite[Section 3]{FPan}.}} on $\MM_{2,n}$, and these are always going to minimal. There are also classes coming from ``endoscopy''. }
\end{rem}

\vspace{+16 pt}
\noindent Departement Mathematik \\
\noindent ETH Z\"urich \\
\noindent rahul@math.ethz.ch  

\vspace{+16 pt}
\noindent CNRS and Versailles St-Quentin University \\
\noindent Laboratoire Math\'e{}matique de Versailles \\
\noindent dimitri.zvonkine@uvsq.fr

\vspace{+16 pt}
\noindent Matematiska institutionen \\
\noindent Stockholms universitet \\
\noindent dan.petersen@math.su.se

\end{document}